\documentclass[twoside,11pt]{article}
\usepackage{amsmath,amssymb,amsfonts}
\begin{document}
\setlength{\unitlength}{10mm}
\newcommand{\f}{\frac}
\newtheorem{theorem}{Theorem}[section]
\newtheorem{lemma}{Lemma}[section]
\newtheorem{fact}{Fact}[section]
\newtheorem{problem}{Problem}[section]
\newtheorem{corollary}{Corollary}[section]
\newcommand{\sta}{\stackrel}
\def\GL{\operatorname{GL}}
\def\SL{\operatorname{SL}}
\def\SK{\operatorname{SK}}
\def\GU{\operatorname{GU}}
\def\SU{\operatorname{SU}}
\def\EU{\operatorname{EU}}
\def\Spin{\operatorname{Spin}}
\def\Epin{\operatorname{Epin}}
\def\Sp{\operatorname{Sp}}
\def\Ep{\operatorname{Ep}}
\def\SO{\operatorname{SO}}
\def\EO{\operatorname{EO}}
\def\FU{\operatorname{FU}}
\def\St{\operatorname{St}}
\def\Max{\operatorname{Max}}
\def\Spec{\operatorname{Spec}}
\def\Hom{\operatorname{Hom}}
\def\rk{\operatorname{rk}}
\def\sr{\operatorname{sr}}
\def\asr{\operatorname{asr}}
\def\sic{\operatorname{sc}}
\def\ad{\operatorname{ad}}
\def\id{\operatorname{id}}
\def\A{\operatorname{A}}
\def\B{\operatorname{B}}
\def\C{\operatorname{C}}
\def\D{\operatorname{D}}
\def\E{\operatorname{E}}
\def\F{\operatorname{F}}
\def\G{\operatorname{G}}
\def\K{\operatorname{K}}
\def\jdim{\operatorname{j-dim}}
\def\Co{{\mathbb C}}
\def\Re{{\mathbb R}}
\def\Int{{\mathbb Z}}
\def\Rat{{\mathbb Q}}
\def\Nat{{\mathbb N}}
\def\e{\varepsilon}
\def\map{\longrightarrow}
\def\GF#1{{\mathbb F}_{\!#1}}
\def\p{{\mathfrak p}}
\def\ma{{\mathfrak m}}


\title{Commutator width in Chevalley groups}
\author{Roozbeh Hazrat
\\
{\scriptsize University of Western Sydney}
\\
Alexei Stepanov
\\
{\scriptsize Saint Petersburg State University,}
\\
{\scriptsize Abdus Salam School of Mathematical Sciences, Lahore}
\\
Nikolai Vavilov
\\
{\scriptsize Saint Petersburg State University}
\\
Zuhong Zhang
\\
{\scriptsize Beijing Institute of Technology} }

\date{}
\maketitle
\begin{picture}(5,2)(3,-9)
 \put(0,.5){{\scriptsize Advances in Group Theory
and Applications}}
 \put(0,0){{\scriptsize Porto Cesareo, June 2011}}
\end{picture}
\vspace{-2.75cm}\begin{abstract}
\noindent
The present paper is the [slightly expanded] text of our
talk at the Conference ``Advances in Group Theory and
Applications'' at Porto Cesareo in June 2011.
Our main results assert that [elementary] Chevalley groups
very rarely have finite commutator width. The reason is that
they have very few commutators, in fact, commutators have
finite width in elementary generators. We discuss also
the background, bounded elementary generation, methods of
proof, relative analogues of these results, some positive
results, and possible generalisations.

\vspace{.5cm}\leftline{{AMS subject Classification 2010:
Primary: 20XX; Secondary: 20XX.}}

\vspace{.5cm}\leftline{Keywords and phrases: {\small Chevalley
groups, elementary subgroups, elementary generators, }}
\leftline{{\small commutator
width, relative groups, bounded
generation, standard commutator formulas,}}
\leftline{{\small unitriangular
factorisations}}
\end{abstract}


\section{Introduction}

\par
In the present note we concentrate on the recent results on the
commutator width of Chevalley groups, the width of commutators in
elementary generators, and the corresponding relative results.
In fact, localisation methods used in the proof of these results
have many further applications, both actual and potential: relative
commutator formulas, multiple commutator formulas, nilpotency
of $K_1$, description of subnormal subgroups, description of
various classes of overgroups, connection with excision kernels,
etc. We refer to our surveys \cite{RN,yoga,yoga2} and
to our papers \cite{RH,RN1,BRN,RZ11,RNZ1,RNZ2,RZ12,step,RNZ3,
HSVZmult,HSVZuni} for these and further applications and many
further related references.


\section{Preliminaries}


\subsection{Length and width}

Let $G$ be a group and $X$ be a set of its generators. Usually
one considers symmetric sets, for which $X^{-1}=X$.
\begin{itemize}\itemsep=0pt
\item
The {\bf length\/} $l_X(g)$ of an element $g\in G$ with
respect to $X$ is the minimal $k$ such that $g$ can be expressed
as the product $g=x_1\ldots x_k$, $x_i\in X$.
\item
The {\bf width} $w_X(G)$ of $G$ with respect to $X$ is
the supremum of $l_X(g)$ over all $g\in G$. In the case when
$w_X(G)=\infty$, one says that $G$ does not have bounded word
length with respect to $X$.
\end{itemize}

The problem of calculating or estimating $w_X(G)$ has
attracted a lot of attention, especially when $G$ is one of
the classical-like groups over skew-fields. There are
{\it hundreds\/} of papers which address this
problem in the case when $X$ is either
\begin{itemize}\itemsep=0pt
\item the set of elementary transvections
\item the set of all transvections or ESD-transvections,
\item the set of all unipotents,
\item the set of all reflections or pseudo-reflections,
\item other sets of small-dimensional transformations,
\item a class of matrices determined by their eigenvalues,
such as the set of all involutions,
\item a non-central conjugacy class,
\item the set of all commutators,
\end{itemize}
\noindent
etc., etc. Many further exotic generating sets have been considered,
such as matrices distinct from the identity matrix in one column,
symmetric matrices, etc., etc., etc. We do not make any attempt to
list all such papers, there are simply far too many, and vast
majority of them produce sharp bounds for classes of rings, which
are trivial from our prospective, such as fields, or semi-local
rings.


\subsection{Chevalley groups}

Let us fix basic notation. This notation is explained in
\cite{abe,AS,matsumoto,stein2,stein,abe88,abe89,NV91,vavplot,NV00},
where one can also find many further references.
\begin{itemize}\itemsep=0pt
\item
$\Phi$ is a reduced irreducible root system;
\item
Fix an order on $\Phi$, let $\Phi^+$, $\Phi^-$
and $\Pi=\{\alpha_1,\ldots,\alpha_l\}$ are the sets of positive,
negative and fundamental roots, respectively.
\item
Let $Q(\Phi)$ be the root lattice of $\Phi$,
$P(\Phi)$ be the weight lattice of $\Phi$ and $P$
be any lattice such that $Q(\Phi)\le P\le P(\Phi)$;
\item
$R$ is a commutative ring with 1;
\item
$G=G_P(\Phi,R)$ is the Chevalley group of type
$(\Phi,P)$ over $R$;
\item
In most cases $P$ does not play essential role and
we simply write $G=G(\Phi,R)$ for any Chevalley group of type
$\Phi$ over $R$;
\item
However, when the answer depends on $P$ we usually write
$G_{\sic}(\Phi,R)$ for the simply connected group, for which
$P=P(\Phi)$ and $G_{\ad}(\Phi,R)$ for the adjoint group,
for which $P=Q(\Phi)$;
\item
$T=T(\Phi,R)$ is a split maximal torus of $G$;
\item
$x_{\alpha}(\xi)$, where $\alpha\in\Phi$, $\xi\in R$,
denote root unipotents $G$ elementary with respect to $T$;
\item
$E(\Phi,R)$ is the [absolute] elementary subgroup
of $G(\Phi,R)$, generated by all root unipotents $x_{\alpha}(\xi)$,
$\alpha\in\Phi$, $\xi\in R$;
\item
$E^L(\Phi,R)$ is the subset (not a subgroup!) of
$E(\Phi,R)$, consisting of products of $\le L$ root unipotents
$x_{\alpha}(\xi)$, $\alpha\in\Phi$, $\xi\in R$;
\item
$H=H(\Phi,R)=T(\Phi,R)\cap E(\Phi,R)$ is the
elementary part of the split maximal torus;
\item
$U^{\pm}(\Phi,R)$ is the unipotent radical of
the standard Borel subgroup $B(\Phi,R)$ or its opposite
$B^-(\Phi,R)$. By definition
\begin{align*}
&U(\Phi,R)=\big\langle
x_{\alpha}(\xi),\ \alpha\in\Phi^+,\ \xi\in R
\big\rangle.\\
&U^-(\Phi,R)=\big\langle
x_{\alpha}(\xi),\ \alpha\in\Phi^-,\ \xi\in R
\big\rangle.\\
\end{align*}


\subsection{Chevalley groups versus elementary subgroups}

Many authors not familiar with algebraic groups or algebraic
$K$-theory do not distinguish Chevalley groups and their elementary
subgroups. Actually, these groups are defined dually.
\item
Chevalley groups $G(\Phi,R)$ are [the groups of
$R$-points of] algebraic groups. In other words, $G(\Phi,R)$
is defined as
$$ G(\Phi,R)=\Hom_{\Int}(\Int[G],R), $$
\noindent
where $\Int[G]$ is the affine algebra of $G$. By definition
$G(\Phi,R)$ consists of solutions in $R$ of certain algebraic
equations.
\item
As opposed to that, {\it elementary\/} Chevalley
groups $E(\Phi,R)$ are generated by elementary generators
$$ E(\Phi,R)=
\big\langle x_{\alpha}(\xi),\ \alpha\in\Phi,\ \xi\in R\big\rangle. $$
\noindent
When $R=K$ is a field, one knows relations among these elementary
generators, so that $E(\Phi,R)$ can be defined by generators and
relations. However, in general, the elementary generators are
described by their action in certain representations.
\end{itemize}

By the very construction of these groups $E(\Phi,R)\le G(\Phi,R)$
but, as we shall see, in general $E(\Phi,R)$ can be strictly smaller
than $G(\Phi,R)$ even for fields. The following two facts might
explain, why some authors confuse $E(\Phi,R)$ and $G(\Phi,R)$:
\begin{itemize}\itemsep=0pt
\item
Let $R=K$ be {\it any\/} field. Then
$G_{\sic}(\Phi,K)=E_{\sic}(\Phi,K)$.
\item
Let $R=K$ be an {\it algebraically closed\/} field. Then
$G_{\ad}(\Phi,K)=E_{\ad}(\Phi,K)$.
\end{itemize}

\noindent
However, for a field $K$ that is not algebraically closed one usually
has strict inclusion $E_{\ad}(\Phi,K)<G_{\ad}(\Phi,K)$. Also, as we
shall see, even for principal ideal domains $E_{\sic}(\Phi,R)<G_{\sic}(\Phi,R)$,
in general.


\subsection{Elementary generators}

By the very construction Chevalley groups occur as subgroups
of the general linear group $\GL(n,R)$. Let $e$ be the identity
matrix and $e_{ij}$, $1\le i,j\le n$, be a matrix unit, which
has 1 in position $(i,j)$ and zeros elsewhere. Below we list
what the elementary root unipotents, also known as elementary
generators, look like for classical groups.

\begin{itemize}\itemsep=0pt
\item
In the case $\Phi=\A_l$ one has $n=l+1$. Root unipotents
of $\SL(n,R)$ are [elementary] transvections
$$ t_{ij}(\xi)=e+\xi e_{ij},\qquad 1\le i\neq j\le n,\quad \xi\in R. $$
\item
In the case $\Phi=\D_l$ one has $n=2l$. We number rows and columns
of matrices from $\GL(n,R)$ as follows: $1,\ldots,l,-l,\ldots,-1$.
Then root unipotents of $\SO(2l,R)$ are [elementary] orthogonal
transvections
$$ T_{ij}(\xi)=e+\xi e_{ij}-\xi e_{-j,-i},\qquad
1\le i,j\le -1,\ i\neq\pm j,\quad \xi\in R. $$
\item
In the case $\Phi=\C_l$ also $n=2l$ and we use the same
numbering
of rows and columns as in the even orthogonal case. Moreover,
we denote $\e_i$ the sign of $i$, which is equal to $+1$ for
$i=1,\ldots,l$ and to $-1$ for $i=-1,\ldots,-1$. In $\C_l$
there are two root lengths. Accordingly, root unipotents of
$\Sp(2l,R)$ come in two stocks. Long root unipotents are the
usual linear transvections $t_{i,-i}(\xi)$, $1\le i\le-1$,
$\xi\in R$, while short root unipotents are [elementary]
symplectic transvections
$$ T_{ij}(\xi)=e+\xi e_{ij}-\e_i\e_j\xi e_{-j,-i},\qquad
1\le i,j\le -1,\ i\neq\pm j,\quad \xi\in R. $$
\item
$\bullet$
Finally, for $\Phi=\B_l$ one has $n=2l+1$ and we number rows and
columns of matrices from $\GL(n,R)$ as follows:
$1,\ldots,l,0,-l,\ldots,-1$. Here too there are two root
lengths. The long root elements of the odd orthogonal group
$\SO(2l+1,R)$ are precisely the root elements of the even
orthogonal groups, $T_{ij}(\xi)$, $i\neq \pm j$, $i,j\neq 0$,
$\xi\in R$. The short root elements have the form
$$ T_{i0}(\xi)=e+\xi e_{i0}-2\xi e_{-i,0}-\xi^2 e_{i,-1},
\qquad i\neq 0,\quad \xi\in R. $$
\end{itemize}

It would be only marginally more complicated to specify
root elements of spin groups and exceptional groups, in their
minimal faithful representations, see \cite{NV00,NV08}.


\subsection{Classical cases}

Actually, most of our results are already new for classical groups.
Recall identification of Chevalley groups and elementary Chevalley
groups for the classical cases. The second column of the following
table lists traditional notation of classical groups, according
to types: $\A_l$ the special linear group, $\B_l$ the odd orthogonal
group, $\C_l$ the symplectic group, and $\D_l$ the even orthogonal
group. These groups are defined by algebraic equations. Orthogonal
groups are not simply connected, the corresponding simply connected
groups are the spin groups. The last column lists the names their
elementary subgroups, generated by the elementary generators listed
in the preceding subsection.

$$ \begin{matrix}
\Phi\qquad\qquad&G(\Phi,R)\hfill&E(\Phi,R)\hfill\\
\noalign{\vskip 10truept}
\A_l\qquad\hfill&\SL(l+1,R)\hfill&E(l+1,R)\hfill\\
\noalign{\vskip 10truept}
\B_l\qquad\hfill&\Spin(2l+1,R)\qquad\hfill&\Epin(2l+1,R)\hfill\\
\noalign{\vskip 5truept}
&\SO(2l+1,R)\hfill&\EO(2l+1,R)\hfill\\
\noalign{\vskip 10truept}
\C_l\qquad\hfill&\Sp(2l,R)\hfill&\Ep(2l,R)\hfill\\
\noalign{\vskip 10truept}
\D_l\qquad\hfill&\Spin(2l,R)\hfill&\Epin(2l,R)\hfill\\
\noalign{\vskip 5truept}
&\SO(2l,R)\hfill&\EO(2l,R)\hfill\\
\end{matrix} $$

Orthogonal groups [and spin groups] in this table are the
{\it split\/} orthogonal groups. {\it Split\/} means that
they preserve a bilinear/quadratic form of maximal Witt index.
In the case of a field the group $\EO(n,K)$ was traditionally
denoted by $\Omega(n,K)$ and called the kernel of spinor
norm. Since the group $\SO(n,K)$ is not simply connected,
in general $\Omega(n,K)$ is a proper subgroup of $\SO(n,K)$.


\subsection{Dimension of a ring}

Usually, dimension of a ring $R$ is defined as the length $d$
of the longest strictly ascending chain of ideals $I_0<I_1<\ldots<I_d$
of a certain class.

\begin{itemize}\itemsep=0pt
\item
The most widely known one is the Krull dimension $\dim(R)$
defined in terms of chains of prime ideals of $R$. Dually, it
can be defined as the combinatorial dimension of $\Spec(R)$,
considered as a topological space with Zariski topology.
\end{itemize}

Recall, that the combinatorial dimension $\dim(X)$ of a topological
space $X$ is the length of the longest {\it descending\/}
chain of its {\it irreducible\/} subspaces $X_0>X_1,\ldots>X_d$.
Thus, by definition, $$\dim(R)=\dim(\Spec(R)).$$
\par\smallskip
However, we mostly use the following more accurate notions of
dimension.

\begin{itemize}\itemsep=0pt
\item
The Jacobson dimension $\jdim(R)$ of $R$ is defined in
terms of $j$-ideals, in other words, those prime ideals, which are
intersections of maximal ideals. Clearly, $\jdim(R)$ coincides
with the combinatorial dimension of the {\it maximal\/} spectrum
of the ring $R$, by definition, $\jdim(R)=\dim(\Max(R))$
\end{itemize}

\noindent
Define dimension $\delta(X)$ of a topological space $X$
as the smallest integer $d$ such that $X$ can be expressed as a {\it
finite\/} union of Noetherian topological spaces of dimension
$\le d$. The trick is that these spaces do not have to be closed
subsets of $X$.

\begin{itemize}\itemsep=0pt
\item
The Bass---Serre dimension of a ring $R$ is defined as
the dimension of its maximal spectrum, $\delta(R)=\delta(\Max(R))$.
\end{itemize}

Bass---Serre dimension has many nice properties, which make it
better adapted to the study of problems we consider.
First, a ring is semilocal iff $\delta(R)=0$ (recall that
a commutative ring $R$ is called semilocal if it has finitely
many maximal ideals). Another usefull property is Serre's theorem asserting
that $\delta(R[t])=\dim(R)+1$.


\subsection{Stability conditions}

Mostly, stability conditions are defined in terms of stability of
rows, or columns. In this note we only refer to Bass' stable rank,
first defined in \cite{bass64}. We will denote the [left]
$R$-module of rows of length $n$ by ${}^n\!R$, to distinguish
it from the [right] $R$-module $R^n$ of columns of height $n$.
\par
A row $(a_1,\ldots,a_n)\in {}^{n}\!R$ is called {\it unimodular\/},
if its components $a_1,\ldots,a_n$ generate $R$ as a right ideal,
$$ a_1R+\ldots+a_nR=R. $$
\noindent
or, what is the same, if there exist such
$b_1,\ldots,b_n\in R$ that
$$ a_1b_1+\ldots+a_nb_n=1. $$
\par
The {\it stable rank\/} $\sr(R)$ of the ring $R$ is the smallest
such $n$ that every unimodular row $(a_1,\ldots,a_{n+1})$ of
length $n+1$ is {\it stable\/}. In other words, there exist
elements $b_1,\ldots b_n\in R$ such that the row
$$ (a_1+a_{n+1}b_1,a_2+a_{n+1}b_2,\ldots,a_n+a_{n+1}b_n) $$
\noindent
of length $n$ is unimodular. If no such $n$ exists, one writes
$\sr(R)=\infty$.
\par
In fact, stable rank is a more precise notion of dimension of
a ring, based on linear algebra, rather than chains of ideals.
It is shifted by 1 with respect to the classical notions of
dimension. The basic estimate of stable rank is Bass' theorem,
asserting that $\sr(R)\le\delta(R)+1$.
\par
Especially important in the sequel is the condition $\sr(R)=1$.
A ring $R$ has stable rank 1 if for any $x,y\in R$ such that
$xR+yR=R$ there exists a $z\in R$ such that $(x+yz)R=R$. In fact,
rings of stable rank 1 are weakly finite (one-sided inverses
are automatically two-sided), so that this last condition is
equivalent to invertibility of $x+yz$. Rings of stable rank
1 should be considered as a class of 0-dimensional rings, in
particular, all semilocal rings have stable rank 1. See
\cite{vaser84} for many further examples and references.


\subsection{Localisation}

Let, as usual, $R$ be a commutative ring with
1, $S$ be a multiplica\-tive system in $R$ and $S^{-1}R$ be the
corresponding localisation. We will mostly use localisation
with respect to the two following types of multiplicative
systems.

\begin{itemize}\itemsep=0pt
\item
Principal localisation: the multiplicative system $S$
is generated by a non-nilpotent element $s\in R$, viz.
$S=\langle s\rangle=\{1,s,s^2,\ldots\}$. In this case we usually
write $\langle s\rangle^{-1}R=R_s$.
\item
Maximal localisation: the multiplicative system $S$
equals $S=R\setminus\ma$, where $\ma\in\Max(R)$ is a maximal
ideal in $R$. In this case we usually write
$(R\setminus\ma)^{-1}R=R_\ma$.
\end{itemize}

We denote by $F_S:R\map S^{-1}R$ the canonical ring homomorphism
called the localisation homomorphism. For the two special cases
mentioned above, we write $F_s:R\map R_s$ and $F_\ma:R\map R_\ma$,
respectively.
\par
Both $G(\Phi,\underline{\ \ })$ and $E(\Phi,\underline{\ \ })$
commute with direct limits. In other words, if $R=\varinjlim R_i$,
where $\{R_i\}_{i\in I}$ is an inductive system of rings, then
$G(\Phi,\varinjlim R_i)=\varinjlim G(\Phi,R_i)$ and the same
holds for $E(\Phi,R)$. Our proofs crucially depend on this
property, which is mostly used in the two following situations.

\begin{itemize}\itemsep=0pt
\item
First, let $R_i$ be the inductive system of all finitely
generated subrings of $R$ with respect to inclusion. Then
$X=\varinjlim X(\Phi,R_i)$, which reduces most of the proofs
to the case of Noetherian rings.
\item
Second, let $S$ be a multiplicative system in $R$ and
$R_s$, $s\in S$, the inductive system with respect to the localisation
homomorphisms: $F_{t}:R_s\map R_{st}$. Then
$X(\Phi,S^{-1}R)=\varinjlim X(\Phi,R_s)$, which allows to reduce
localisation with respect to any multiplicative system to principal
localisations.
\end{itemize}


\subsection{$\K_1$-functor}

The starting point of the theory we consider is the following result,
first obtained by Andrei Suslin \cite{Sus} for $\SL(n,R)$, by Vyacheslav
Kopeiko \cite{kopeiko} for symplectic groups, by Suslin and Kopeiko
\cite{suskop} for even orthogonal groups and by Giovanni Taddei
\cite{taddei} in general.

{\theorem Let\/ $\Phi$ be a reduced irreducible root system such
that\/ $\rk(\Phi)\ge 2$. Then for any commutative ring\/ $R$ one
has\/ $E(\Phi,R)\trianglelefteq G(\Phi,R)$.}

\bigskip

In particular, the quotient
$$ K_1(\Phi,R)=G_{\sic}(\Phi,R)/E_{\sic}(\Phi,R) $$
\noindent
is not just a pointed set, it is a group. It is called
$K_1$-functor.
\par
The groups $G(\Phi,R)$ and $E(\Phi,R)$ behave functorially with
respect to both $R$ and $\Phi$. In particular, to an embedding
of root systems $\Delta\subseteq\Phi$ there corresponds the
map $\phi: G(\Delta,R)\map G(\Phi,R)$ of the corresponding
[simply connected] groups, such that
$\phi(E(\Delta,R))\le G(\Phi,R)$. By homomorphism theorem
it defines the stability map $\phi:K_1(\Delta,R)\map K_1(\Phi,R)$.
\par
In the case $\Phi=\A_l$ this $K_1$-functor specialises to the
functor 
$$ \SK_1(n,R)=\SL(n,R)/E(n,R), $$
\noindent
rather than the usual linear $K_1$-functor 
$K_1(n,R)=\GL(n,R)/E(n,R)$. In examples below we also mention
the corresponding stable $K_1$-functors, which are defined as
limits of $K_1(n,R)$ and $\SK_1(n,R)$ under stability embeddings,
as $n$ tends to infinity:
$$ \SK_1(R)=\varinjlim\SK_1(n,R),\qquad 
K_1(R)=\varinjlim K_1(n,R). $$
\par
Another basic tool are stability theorems, which assert that
under some assumptions on $\Delta,\Phi$ and $R$ stability maps
are surjective or/and injective. We do not try to precisely
state stability theorems for Chevalley groups, since they
depend on various analogues and higher versions of stable rank,
see in particular \cite{stein,plotkin2,plotkin1,plotkin}.
\par
However, to give some feel, we state two classical results
pertaining to the case of $\SL(n,R)$. These results, which
are due to Bass and Bass---Vaserstein, respectively, are known
as surjective stability of $K_1$ and injective stability of $K_1$.
In many cases they allow to reduce problems about groups of
higher ranks, to similar problems for groups of smaller rank.

{\theorem For any $n\ge\sr(R)$ the stability map
$$K_1(n,R)\map K_1(n+1,R)$$ is surjective. In other words,
$$ \SL(n+1,R)=\SL(n,R) E(n+1,R). $$
\vskip -2truecm}

{\theorem For any $n\ge\sr(R)+1$ the stability map
$$K_1(n,R)\map K_1(n+1,R)$$ is injective. In other words,
$$ \SL(n,R)\cap E(n+1,R)=E(n,R). $$
\vskip -2truecm}


\subsection{$\K_1$-functor: trivial or non-trivial}

Usually, $K_1$-functor is non-trivial. But in some important
cases it is trivial. Let us start with some obvious examples.

\begin{itemize}\itemsep=0pt
\item
$R=K$ is a field.
\item
More generally, $R$ is semilocal
\item
$R$ is Euclidean
\item
It is much less obvious that $K_1$ does not have to be trivial
even for principal ideal rings. Let us cite two easy examples
discovered by Ischebeck \cite{ische} and by Grayson and
Lenstra \cite{Grayson}, respectively.
\item
Let $K$ be a field of algebraic functions of one
variable with a perfect field of constants $k$.
Then the ring $R=K\otimes_k k(x_1,\ldots,x_m)$ is a principal
ideal ring. If, moreover, $m\ge 2$, and the genus of $K$
is distinct from $0$, then $\SK_1(R)\neq 1$.
\item
Let $R=\Int[x]$, and $S\subseteq R$ be the multiplicative
subsystem of $R$ generated by cyclotomic polynomials $\Phi_n$,
$n\in\Nat$. Then $S^{-1}R$ is a principal ideal ring such that
$\SK_1(S^{-1}R)\neq 1$.
\end{itemize}

This is precisely why over a Euclidean ring it is somewhat easier
to find Smith form of a matrix, than over a principal ideal ring.
\par
However, there are some further examples, when $K_1$ is trivial.
Usually, they are very deep. The first example below is part of
the [almost] positive solution of the congruence subgroup problem
by Bass---Milnor---Serre and Matsumoto \cite{bassmilnorserre,
matsumoto}. The second one is the solution of $K_1$-analogue of
Serre's problem by Suslin \cite{Sus}.

\begin{itemize}\itemsep=0pt
\item
$R={\cal O}_S$ is a Hasse domain.
\item
$R=K[x_1,\ldots,x_m]$ is a polynomial ring over a field.
\end{itemize}


\subsection{$\K_1$-functor, abelian or non-abelian}

Actually, $K_1(\Phi,R)$ is not only non-trivial. Oftentimes, it is
even non-abelian. The first such examples were constructed by
Wilberd van der Kallen \cite{kallen89} and Anthony Bak \cite{B4}.
In both cases proofs are of topological nature and use homotopy
theory.

\begin{itemize}\itemsep=0pt
\item
Wilberd van der Kallen \cite{kallen89} constructs a
number of examples, where $K_1(n,R)$ is non-abelian. For instance,
$$ R=\Re[x_1,x_2,y_1,y_2,y_3,y_4]/
(x_1^2+x_2^2=y_1^2+y_2^2+y_3^2+y_4^2=1) $$
\noindent
is a 4-dimensional ring for which $[\SL(4,R),\SL(4,R)]\not\le E(4,R)$.
In fact, in this case even
$$ [\SL(2,R),\SL(4,R)]\not\le\GL(3,R) E(4,R). $$
\item
Anthony Bak \cite{B4} constructs examples of [finite
dimensional] subrings $R$ in the rings of continuous functions
$\Re^X$ and $\Co^X$ on certain topological spaces $X$, for
which not only $K_1(n,R)$, $n\ge 3$, is non-abelian, but even
its nilpotency class can be arbitrarily large.
\end{itemize}

The question arises, as to how non-abelian $K_1(\Phi,R)$ may be.
For finite dimensional rings this question was answered by
Anthony Bak \cite{B4} for $\SL(n,R)$, for other even classical
groups by the first author \cite{RH} and for all Chevalley groups
by the first and the third authors \cite{RN1}.

{\theorem Let $\Phi$ be a reduced irreducible root system such that $\rk(\Phi)\ge 2$. Further let $R$ be a commutative ring
of Bass---Serre dimension $\delta(R)=d<\infty$. Then $K_1(\Phi,R)$
is nilpotent of class $\le d+1$.}

\bigskip

This theorem relies on a version of localisation method which
Bak called localisation-completion \cite{B4}. This method turned
out to be crucial for the proof of results we discuss in the
present paper, see \cite{RN,yoga} for more historical background and
an introduction to this method in non-technical terms.


\section{Main problems}


\subsection{Statement of the main problems}

In this paper we discuss the following problem.

{\problem Estimate the width of $E(\Phi,R)$ with respect to
the set of elementary commutators
$$ X=\big\{[x,y]=xyx^{-1}y^{-1},
\ x\in G(\Phi,R),\ y\in E(\Phi,R)\big\}. $$}

Observe, that one could not have taken the set
$$ X=\big\{[x,y]=xyx^{-1}y^{-1},\ x,y\in G(\Phi,R)\big\} $$
\noindent
here, since $K_1(\Phi,R)$ maybe non-abelian.

It turns out that this problem is closely related to the
following problem.

{\problem Estimate the width of $E(\Phi,R)$ with respect to
the set of elementary generators
$$ X=\big\{x_{\alpha}(\xi),\ \alpha\in\Phi,\ \xi\in R\big\}. $$}

The answer in general will be highly unexpected, so we start
with discussion of classical situations.


\subsection{The group $\SL(2,R)$}

Let us mention one assumption that is essential in what follows.

When $R$ is Euclidean, expressions of matrices in $\SL(2,R)$
as products of elementary transvections correspond to
continued fractions. Division chains in $\Int$ can be arbitrarily
long, it is classically known that two consecutive Fibonacci
numbers provide such an example. Thus, we get.

{\fact $\SL(2,\Int)$ does not have bounded length with respect
to the elementary generators.}

\bigskip

Actually, behavious of the group $\SL(2,R)$ is exceptional in
more than one respect. Thus, the groups $E(n,R)$, $n\ge 3$ are
perfect. The group $E(2,R)$ is usually not.

{\fact $[\SL(2,\Int),\SL(2,\Int)]$ has index $12$ in $\SL(2,\Int)$.}

\begin{itemize}\itemsep=0pt
\item
Thus, in the sequel we always assume that $\rk(\Phi)\ge 2$.
\item
In fact, it is material for most of our results that the group
$E(\Phi,R)$ is perfect. It usually is, the only counter-examples
in rank $\ge 2$ stemming from the fact that $\Sp(4,\GF{2})$ and
$G(\G_2,\GF{2})$ are not perfect. Thus, in most cases one should
add proviso that $E(\Phi,R)$ is actually perfect, which amounts
to saying that $R$ does not have residue field $\GF{2}$ for
$\Phi=\B_2,\G_2$.
\end{itemize}
\noindent
The reader may take these two points as standing assumptions for
the rest of the note.


\subsection{The answers for fields}

The following result easily follows from Bruhat decomposition.

{\theorem The width of $G_{\sic}(\Phi,K)$ with respect to
the set of elementary generators is $\le 2|\Phi^+|+4\rk(\Phi)$.}

\bigskip

Rimhak Ree \cite{ree} observed that the commutator width of
semisimple algebraic groups over an algebraically closed fields
equals 1. For fields containing $\ge 8$ elements the following
theorems were established by Erich Ellers and Nikolai Gordeev
\cite{EG} using Gauss decomposition with prescribed semi-simple
part \cite{CEG}. On the other hand, for very small fields these
theorems were recently proven by Martin Liebeck, Eamonn O'Brien,
Aner Shalev, and Pham Huu Tiep \cite{LOST10,LOST11}, using
explicit information about maximal subgroups and very delicate
character estimates.

Actually, the first of these theorems in particular
completes the answer to Ore conjecture, whether any element
of a [non-abelian] finite simple group is a single
commutator.

{\theorem The width of $E_{\ad}(\Phi,K)$ with respect to
commutators is $1$.}
{\theorem The width of $G_{\sic}(\Phi,K)$ with respect to
commutators is $\le 2$.}


\subsection{The answers for semilocal rings}

The following results were recently published by by Andrei
Smolensky, Sury and the third author \cite{SSV,VSS}. Actually,
their proofs are easy combinations of Bass' surjective stability
\cite{bass64} and Tavgen's rank reduction theorem \cite{tavgen90}.
The second of these decompositions, the celebrated Gauss
decomposition, was known for semilocal rings, the first one was
known for $\SL(n,R)$, see \cite{DV89}, but not in general.

{\theorem Let $\sr(R)=1$. Then the
$$ E(\Phi,R)=U^+(\Phi,R)U^-(\Phi,R)U^+(\Phi,R)U^-(\Phi,R). $$
\vskip -2truecm}

{\corollary Let $\sr(R)=1$. Then the width of $E(\Phi,R)$
with respect to the set of elementary generators is at most
$M=4|\Phi^+|$.}

{\theorem Let $\sr(R)=1$. Then the
$$ E(\Phi,R)=U^+(\Phi,R)U^-(\Phi,R)H(\Phi,R)U(\Phi,R). $$
\vskip -2truecm}

{\corollary Let $\sr(R)=1$. Then the width of $E(\Phi,R)$
with respect to the set of elementary generators is at most
$M=3|\Phi^+|+4\rk(\Phi)$.}

\bigskip

In particular, the width of $E(\Phi,R)$ over a ring with
$\sr(R)=1$ with respect to commutators is always bounded,
but its explicit calculation is a non-trivial task.
Let us limit ourselves with the following result by Leonid
Vaserstein and Ethel Wheland \cite{VW1,VW2}.

{\theorem Let $\sr(R)=1$. Then the width of $E(n,R)$, $n\ge 3$,
with respect to commutators is $\le 2$.}

\bigskip

There are also similar results by You Hong, Frank Arlinghaus
and Leonid Vaserstein \cite{you94,AVY} for other classical
groups, but they usually assert that the commutator width
is $\le 3$ or $\le 4$, and sometimes impose stronger stability
conditions such as $\asr(R)=1$, $\Lambda\sr(R)=1$, etc.
\par
The works by Nikolai Gordeev and You Hong,
where similar results are established for exceptional groups
over local rings [subject to some mild restrictions on their
residue fields] are still not published.


\subsection{Bounded generation}

Another nice class of rings, for which one may expect positive
answers to the above problems, are Dedekind rings of arithmetic
type.
\par
Let $K$ be a global field, i.e.\ either a finite algebraic extension
of $\Rat$, or a field of algebraic functions of one variable over
a finite field of constants, and further let $S$ be a finite set of
(non-equivalent) valuations of $K$, which is non-empty in the
functional case, and which contains all Archimedian valuations in the
number case. For a non-Archimedian valuation $\p$ of the field $K$
we denote by $v_{\p}$ the corresponding exponent. As usual,
$R={\mathcal O}_S$ denotes the ring, consisting of $x\in K$ such
that $v_{\p}(x)\ge 0$ for all valuations $\p$ of $K$, which do not
belong to $S$. Such a ring ${\mathcal O}_S$ is known as the Dedekind
ring of arithmetic type, determined by the set of valuations $S$
of the field $K$. Such rings are also called Hasse domains, see,
for instance, \cite{bassmilnorserre}. Sometimes one has to require
that $|S|\ge 2$, or, what is the same, that the multiplicative
group ${\mathcal O}_S^*$ of the ring ${\mathcal O}_S$ is infinite.
\par
Bounded generation of $\SL(n,{\mathcal O}_S)$, $n\ge 3$, was established by
David Carter and Gordon Keller in \cite{CK83,CK84,CK85a,CK85b,CKP85},
see also the survey by Dave Witte Morris \cite{Mo} for a modern
exposition. The general
case was solved by Oleg Tavgen \cite{tavgen90,tavgen92}. The result
by Oleg Tavgen can be stated in the following form due to the
[almost] positive solution of the congruence subgroup problem
\cite{bassmilnorserre,matsumoto}.

{\theorem Let ${\mathcal O}_S$ be a Dedekind ring of arithmetic
type, $\rk(\Phi)\ge 2$. Then the elementary Chevalley group
$G(\Phi,{\mathcal O}_S)$ has bounded length with
respect to the elementary generators.}

\bigskip

In Section~6 we discuss what this implies for the commutator width.

See also the recent works by Edward Hinson \cite{hinson},
Loukanidis and Murty \cite{LM,Mu},
Sury \cite{sury}, Igor Erovenko and Andrei Rapinchuk
\cite{erovenkothesis,ER01,ER06}, for different proofs,
generalisations and many further references, concerning bounded
generation.


\subsection{van der Kallen's counter-example}

However, all hopes for positive answers in general are
completely abolished by the following remarkable result
due to Wilberd van der Kallen \cite{kallen82}.

{\theorem The group $\SL(3,\Co[t])$ does not have bounded word
length with respect to the elementary generators.}

\bigskip

It is an amazing result, since $\Co[t]$ is Euclidean.
Since $\sr(\Co[t])=2$ we get the following corollary

{\corollary None of the groups $\SL(n,\Co[t])$, $n\ge 3$,
has bounded word length with respect to the elementary generators.}

\bigskip

See also \cite{erovenko} for a slightly easier proof of a slightly
stronger result. Later Dennis and Vaserstein \cite{DV89} improved van
der Kallen's result to the following.

{\theorem The group $\SL(3,\Co[t])$ does not have bounded word
length with respect to the commutators.}

\bigskip

Since for $n\ge 3$ every elementary matrix is a commutator,
this is indeed stronger, than the previous theorem.


\section{Absolute commutator width}

Here we establish an amazing relation between Problems 3.1
and 3.2.


\subsection{Commutator width in $\SL(n,R)$}

The following result by Alexander Sivatsky and the second author
\cite{SiSt} was a major breakthrough.

{\theorem Let\/ $n\ge 3$ and let $R$ be a Noetherian ring such
that\/ $\dim\Max(R)=d<\infty$. Then there exists a natural
number\/ $N=N(n,d)$ depending only on\/ $n$ and\/ $d$ such that each
commutator\/ $[x,y]$ of elements\/ $x\in E(n,R)$
and\/ $y\in\SL(n,R)$ is a product of at most\/ $N$ elementary
transvections.}

\bigskip

Actually, from the proof in \cite{SiSt} one can derive an
efficient upper bound on $N$, which is a {\it polynomial\/}
with the leading term $48n^6d$.

It is interesting to observe that it is already non-trivial
to replace here an element of $\SL(n,R)$ by an element of
$\GL(n,R)$. Recall, that a {\it ring of geometric origin\/}
is a localisation of an affine algebra over a field.

{\theorem Let\/ $n\ge 3$ and let $R$ be a ring of geometric
origin. Then there exists a natural number\/ $N$ depending
only on\/ $n$ and\/ $R$ such that each commutator\/ $[x,y]$
of elements\/ $x\in E(n,R)$ and\/ $y\in\GL(n,R)$ is a product
of at most\/ $N$ elementary transvections.}

\bigskip

Let us state another interesting variant of the Theorem 4.1,
which may be considered as its stable version. Its proof
crucially depends on the Suslin---Tulenbaev proof of the
Bass---Vaserstein theorem, see \cite{sustul}.

{\theorem Let\/ $n\ge\sr(R)+1$. Then there exists a natural
number\/ $N$ depending only on\/ $n$ such that each
commutator\/ $[x,y]$ of elements\/ $x,y\in\GL(n,R)$ is a
product of at most\/ $N$ elementary transvections.}

\bigskip


Actually, \cite{SiSt} contains many further interesting results,
such as, for example, analogues for the Steinberg groups
$\St(n,R)$, $n\ge 5$. However, since this result depends on
the centrality of $K_2(n,R)$ at present there is no hope
to generalise it to other groups.


\subsection{Decomposition of unipotents}

The proof of Theorem 4.1 in \cite{SiSt} was based on a combination
of localisation and decomposition of unipotents \cite{SV}.
Essentially, in the simplest form decomposition of unipotents
gives finite polynomial expressions of the conjugates
$$ gx_{\alpha}(\xi)g^{-1},\qquad
\alpha\in\Phi,\ \xi\in R,\ g\in G(\Phi,R), $$
\noindent
as products of factors sitting in proper parabolic subgroups,
and, in the final count, as products of elementary generators.
\par
Roughly speaking, decomposition of unipotents allows to plug
in explicit polynomial formulas as the induction base ---
which is the most difficult part of all localisation proofs! ---
instead of messing around with the length estimates in the
conjugation calculus.
\par
To give some feel of what it is all about, let us state an
immediate corollary of the Theme of \cite{SV}. Actually,
\cite{SV} provides explicit polynomial expressions of the
elementary factors, rather than just the length estimate.

{\fact Let\/ $R$ be a commutative ring and $n\ge 3$. Then
any transvection of the form
$gt_{ij}(\xi)g^{-1}$, $1\le i\neq j\le n$, $\xi\in R$,
$g\in\GL(n,R)$ is a product of at most $4n(n-1)$ elementary
transvections.}

\bigskip

It is instructive to compare this bound with the bound resulting
from Suslin's proof of Suslin's normality theorem \cite{Sus}.
Actually, Suslin's direct factorisation method is more general,
in that it yields elementary factorisations of a broader class
of transvections. On the other hand, it is less precise, both
factorisations coincide for $n=3$, but asymptotically
factorisation in Fact 4.1 is better.

{\fact Let\/ $R$ be a commutative ring and $n\ge 3$. Assume
that $u\in R^n$ is a unimodular column and $v\in{}^n\!R$ be any
row such that $vu=0$. Then the transvection $e+uv$ is a product
of at most $n(n-1)(n+2)$ elementary transvections.}

\bigskip

Let us state a counterpart of the Theorem 4.1 that results from
the Fact 4.1 alone, {\it without\/} the use of localisation.
This estimate works for {\it arbitrary\/} commutative rings,
but depends on the length of the elementary factor. Just wait
until subsection 4.5!

{\theorem Let\/ $n\ge 3$ and let $R$ be a commutative ring. Then
there exists a natural number\/ $N=N(n,M)$ depending only on\/ $n$
and\/ $M$ such that each commutator\/ $[x,y]$ of elements\/
$x\in E^M(n,R)$ and\/ $y\in\SL(n,R)$ is a product of at
most\/ $N$ elementary transvections.}

\bigskip

It suffices to expand a commutator $[x_1\ldots x_M,y]$, where
$x_i$ are elementary transvections, with the help of the
commutator identity $[xz,y]={}^x[z,y]\cdot[x,y]$, and take
the upper bound $4n(n-1)+1$ for each of the resulting
commutators $[x_i,y]$. One thus gets $N\le M^2+4n(n-1)M$.

However, such explicit formulas are only available for linear
and orthogonal groups, and for exceptional groups of types
$\E_6$ and $\E_7$. Let us state the estimate resulting from
the proof of \cite{NV00}, Theorems 4 and 5.

{\fact Let\/ $R$ be a commutative ring and $\Phi=\E_6,\E_7$.
Then any root element of the form $gx_{\alpha}(\xi)g^{-1}$,
$\alpha\in\Phi$, $\xi\in R$, $g\in G(\Phi,R)$ is a product
of at most $4\cdot16\cdot27=1728$ elementary root unipotents
in the case of $\Phi=\E_6$ and of at most $4\cdot27\cdot56=6048$
elementary root unipotents in the case of $\Phi=\E_7$.}

\bigskip

Even for symplectic groups --- not to say
for exceptional groups of types $\E_8,\F_4$ and $\G_2$! ---
it is only known that the elementary groups are generated by
root unipotents of certain classes, which afford reduction to
smaller ranks, but no explicit polynomial factorisations are
known, and even no polynomial length estimates.
\par
This is why generalisation of Theorem 4.1 to Chevalley groups
requires a new idea.


\subsection{Commutator width of Chevalley groups}

Let us state the main result of \cite{SV10}. While the main idea
of proof comes from the work by Alexander Sivatsky and the second
author \cite{SiSt}, most of the actual calculations are refinements
of conjugation calculus and commutator calculus in Chevalley
groups, developed by the first and the third authors in \cite{RN1}.

{\theorem Let\/ $G=G(\Phi,R)$ be a Chevalley group of rank\/
$l\ge 2$ and let\/ $R$ be a ring such that\/ $\dim\Max(R)=d<\infty$.
Then there exists a natural number\/ $N$ depending only on\/ $\Phi$
and\/ $d$ such that each commutator\/ $[x,y]$ of
elements\/ $x\in G(\Phi,R)$ and\/ $y\in E(\Phi,R)$ is a product
of at most\/ $N$ elementary root unipotents.}

\bigskip

Here we cannot use decomposition of unipotents. The idea of the
second author was to use the {\it second localisation\/} instead.
As in \cite{SiSt} the proof starts with the following lemma, where
$M$ has the same meaning as in Subsection 3.4.

{\lemma Let\/ $d=\dim(\Max(R))$ and\/ $x\in\G(\Phi,R)$.
Then there exist\/ $t_0,\dots,t_k\in R$, where\/ $k\le d$,
generating\/ $R$ as an ideal and such that\/
$F_{t_i}(x)\in\E^{M}(\Phi,R_{t_i})$ for all\/ $i=0,\dots,k$.}

\bigskip

Since $t_0,\ldots,t_k$ are unimodular, their powers also are,
so that we can rewrite $y$ as a product of $y_i$, where each
$y_i$ is congruent to $e$ modulo a high power of $t_i$. In the
notation of the next section this means that $y_i\in E(\Phi,R,t_i^mR)$.
\par
When the ring $R$ is Noetherian, $G(\Phi,R,t_i^mR)$ injects into
$\G(\Phi,R_{t_i})$ for some high power $t_i^m$. Thus, it suffices
to show that $F_{t_i}([x,y_i])$ is a product of bounded number
of elementary factors without denominators in $E(\Phi,R_{t_i})$.
This is the first localisation.
\par
The second localisation consists in applying the same argument
again, this time in $R_{t_i}$. Applying Lemma 4.1 once more we
can find $s_0,\ldots,s_d$ forming a unimodular row, such that
the images of $y_i$ in $E(\Phi,R_{t_is_j})$ are products of at
most $M$ elementary root unipotents with denominators $s_j$.
Decomposing $F_{s_j}(x)\in E(\Phi,R_{s_j})$ into a product
of root unipotents, and repeatedly applying commutator
identities, we eventually reduce the proof to proving that
the length of each commutator of the form
$$ \Big[x_{\alpha}\Big(\frac{t_i^l}{s_j}a\Big),
x_{\beta}\Big(\frac{s_j^n}{t_i}b\Big)\Big] $$
\noindent
is bounded.


\subsection{Commutator calculus}

Conjugation calculus and commutator calculus consists in
rewriting conjugates/commuta\-tors with denominators as products
of elementary generators {\it without\/} denominators.

Let us state a typical technical result, the base of induction
of the commutator calculus.

{\lemma Given $s,t\in R$ and $p,q,k,m\in\Nat$, there exist
$l,m\in\Nat$ and $L=L(\Phi)$ such that
$$ \Big[x_{\alpha}\Big(\frac{t^l}{s^k}a\Big),
x_{\beta}\Big(\frac{s^n}{t^m}b\Big)\Big]\in E^L(\Phi,s^pt^qR). $$}

A naive use of the Chevalley commutator formula gives
$L\le 585$ for simply laced systems, $L\le 61882$ for doubly
laced systems and $L\le 797647204$ for $\Phi=\G_2$. And this is
just the first step of the commutator calculus!

Reading the proof sketched in the previous subsection upwards,
and repeatedly using commutator identities, we can eventually
produce bounds for the length of commutators, ridiculous as
they can be.
\par
Recently in \cite{HSVZuni} the authors succeeded in producing
a similar proof for Bak's unitary groups, see
\cite{HOM,knus,BV3,RN} and references there. The situation here
is in many aspects more complicated than for Chevalley groups.
In fact, Bak's unitary groups are not always algebraic, and all
calculations should be inherently carried through in terms of
{\it form ideals\/}, rather then ideals of the ground ring.
Thus, the results of \cite{HSVZuni} heavily depend on the
{\it unitary\/} conjugation calculus and commutator calculus,
as developed in \cite{RH,RNZ1}.


\subsection{Universal localisation}

Now something truly amazing will happen. Some two years ago the
second author noticed that the width of commutators is bounded
by a universal constant that depends on the type of the group
alone, see \cite{step}. Quite remarkably, one can obtain a
length bound that does not depend either on the dimension of
the ring, or on the length of the elementary factor.

{\theorem Let\/ $G=G(\Phi,R)$ be a Chevalley group of rank\/
$l\ge 2$. Then there exists a natural number\/ $N=N(\Phi)$
depending on\/ $\Phi$ alone, such that each commutator\/ $[x,y]$
of elements\/ $x\in G(\Phi,R)$ and\/ $y\in E(\Phi,R)$ is a
product of at most\/ $N$ elementary root unipotents.}

\bigskip

What is remarkable here, is that there is no dependence on $R$
whatsoever. In fact, this bound applies even to infinite
dimensional rings! Morally, it says that in the groups of points
of algebraic groups there are very few commutators.
\par
Here is a very brief explanation of how it works. First of all,
Chevalley groups are representable functors,
$G(\Phi,R)=\Hom(\Int[G],R)$, so that there is a {\it universal
element\/} $g\in\G(\Phi,\Int[G])$, corresponding to
$\id:\Int[G]\map\Int[G]$, which specialises to {\it any\/}
element of the Chevalley group $G(\Phi,R)$ of the same type over
{\it any\/} ring.
\par
But the elementary subgroup $E(\Phi,R)$ is not an algebraic group,
so where can one find universal elements?
\par
The real know-how proposed by the second author consists in
construction of the universal coefficient rings for the
principal congruence subgroups $G(\Phi,R,sR)$ (see the next
section, for the definition), corresponding to the
principal ideals. It turns out that this is enough to
carry through the same scheme of the proof, with bounds that
do not depend on the ring $R$.


\section{Relative commutator width}

In the absolute case the above results on commutator width are
mostly published. In this section we state relative analogues
of these results which are announced here for the first time.


\subsection{Congruence subgroups}

Usually, one defines congruence subgroups as follows. An ideal
$A\trianglelefteq R$ determines the reduction homomorphism
$\rho_{A}:R\map R/A$. Since $G(\Phi,\underline{\ \ })$ is a
functor from rings to groups, this homomorphism induces
reduction homomorphism $\rho_{A}:G(\Phi,R)\map G(\Phi,R/A)$.

\begin{itemize}\itemsep=0pt
\item
The kernel of the reduction homomorphism $\rho_{A}$
modulo $A$ is called the principal congruence subgroup of
level $A$ and is denoted by $G(\Phi,R,A)$.
\item
The full pre-image of the centre of $G(\Phi,R/A)$ with
respect to the reduction homomorphism $\rho_{A}$ modulo $A$
is called the full congruence subgroup of level $A$, and
is denoted by $C(\Phi,R,A)$.
\end{itemize}

But in fact, without assumption that $2\in R^*$ for doubly laced
systems, and without assumption that $6\in R^*$ for $\Phi=\G_2$,
the genuine congruence subgroups should be defined in terms of
admissible pairs of ideals $(A,B)$, introduced by Abe,
\cite{abe,AS,abe88,abe89}, and in terms of form ideals for
symplectic groups. One of these ideals corresponds to short
roots and another one corresponds to long roots.
\par
In \cite{HPV} we introduced a more general notion of congruence
subgroups $G(\Phi,R,A,B)$ and $C(\Phi,R,A,B)$, corresponding to
admissible pairs. Not to overburden the note with technical
details, we mostly tacitly assume that $2\in R^*$ for
$\Phi=\B_l,\C_l,\F_4$ and $6\in R^*$ for $\Phi=\G_2$.
Under these simplifying assumption one has $A=B$ and
$G(\Phi,R,A,B)=G(\Phi,R,A)$ and $C(\Phi,R,A,B)=C(\Phi,R,A)$.
Of course, using admissible pairs/form ideals one can obtained
similar results without any such assumptions.


\subsection{Relative elementary groups}

Let $A$ be an additive subgroup of $R$. Then $E(\Phi,A)$ denotes
the subgroup of $E$ generated by all elementary root unipotents
$x_{\alpha}(\xi)$ where $\alpha\in\Phi$ and $\xi\in A$.  Further,
let $L$ denote a nonnegative integer and let $E^L(\Phi,A)$ denote
the {\it subset\/} of $E(\Phi,A)$ consisting of all products of
$L$ or fewer elementary root unipotents $x_{\alpha}(\xi)$, where
$\alpha\in\Phi$ and $\xi\in A$. In particular, $E^1(\Phi,A)$ is
the set of all $x_{\alpha}(\xi)$, $\alpha\in\Phi$, $\xi\in A$.
\par
In the sequel we are interested in the case where $A=I$ is an
ideal of $R$. In this case we denote by
$$ E(\Phi,R,I)=E(\Phi,I)^{E(\Phi,R)} $$
\noindent
the {\it relative\/} elementary subgroup of level $I$. As a
{\it normal\/} subgroup of $E(\Phi,R)$ it is generated by
$x_{\alpha}(\xi)$, $\alpha\in\Phi$, $\xi\in A$. The following
theorem \cite{stein2,tits,vaser86} lists its generators as a
subgroup.

{\theorem As a subgroup\/ $E(\Phi,R,I)$ is generated by the
elements
$$ z_{\alpha}(\xi,\zeta)=
x_{-\alpha}(\zeta)x_{\alpha}(\xi)x_{-\alpha}(-\zeta), $$
\noindent
where\/ $\xi\in I$ for $\alpha\in\Phi$, while\/ $\zeta\in R$.}

\bigskip

It is natural to regard these generators as the {\it elementary
generators\/} of $E(\Phi,R,I)$. For the special linear group
$\SL(n,{\mathcal O}_S)$, $n\ge 3$, over a Dedekind ring of
arithmetic type Bernhard Liehl \cite{liehl} has proven bounded
generation of the elementary relative subgroups
$E(n,{\mathcal O}_S,I)$ in the generators $z_{ij}(\xi,\zeta)$.
What is remarkable in his result, is that the bound does not
depend on the ideal $I$. Also, he established similar results
for $\SL(2,{\mathcal O}_S)$, provided that ${\mathcal O}_S^*$
is infinite.


\subsection{Standard commutator formula}

The following result was first proven by Giovanni Taddei
\cite{taddei}, Leonid Vaserstein \cite{vaser86} and Eiichi Abe
\cite{abe88,abe89}.

{\theorem Let\/ $\Phi$ be a reduced irreducible root system of
rank\/ $\ge 2$, $R$ be a commutative ring,\/ $I\trianglelefteq R$
be an ideal of $R$. In the case, where\/ $\Phi=\B_2$ or\/ $\Phi=\G_2$
assume moreover that\/ $R$ has no residue fields\/ ${\mathbb F}_{\!2}$
of\/ $2$ elements. Then the following standard commutator formula
holds
$$ \big [G(\Phi,R),E(\Phi,R,I)\big]=
\big [E(\Phi,R),C(\Phi,R,I)\big]=
E(\Phi,R,I). $$}

In fact, in \cite{HPV} we established similar result for
relative groups defined in terms of admissible pairs, rather
then single ideals. Of course, in all cases, except Chevalley
groups of type $\F_4$, it was known before,
\cite{BV3,petrov2,costakeller2}.

With the use of level calculations the following result was
established by You Hong \cite{you92}, by analogy with the
Alec Mason and Wilson Stothers \cite{MAS3,Mason74,MAS1,MAS2}.
Recently the first, third and fourth authors gave another
proof, of this result, in the framework of relative
localisation \cite{RNZ2}, see also \cite{VS8,RZ11,VS10,yoga,
RNZ1,RZ12,yoga2,RNZ3,HSVZmult,step} for many further
analogues and generalisations of such formulas.

{\theorem Let\/ $\Phi$ be a reduced irreducible root system,\/ $rk(\Phi)\ge 2$.
Further, let\/ $R$ be a commutative ring, and\/ $A,B\trianglelefteq R$ be
two ideals of\/ $R$. Then
$$ [E(\Phi,R,A),G(\Phi,R,B)]=[E(\Phi,R,A),E(\Phi,R,B)]. $$}


\subsection{Generation of mixed commutator subgroups}

It is easy to see that the mixed commutator $[E(\Phi,R,A),E(\Phi,R,B)]$
is a subgroup of level $AB$, in other words, it sits between
the relative elementary subgroup $E(\Phi,R,AB)$ and the corresponding
congruence subgroup $G(\Phi,R,AB)$.

{\theorem Let\/ $\Phi$ be a reduced irreducible root system,\/ $rk(\Phi)\ge 2$.
Further, let\/ $R$ be a commutative ring, and\/ $A,B\trianglelefteq R$ be
two ideals of\/ $R$. Then
$$ E(\Phi,R,AB)\le [E(\Phi,R,A),E(\Phi,R,B)]\le
[G(\Phi,R,A),G(\Phi,R,B)]\le G(\Phi,R,AB). $$}

It is not too difficult to construct examples showing that in general
the mixed commutator subgroup
$E(\Phi,R,A),E(\Phi,R,B)]$ can be strictly larger than
$E(\Phi,R,AB)$. The first such examples were constructed by
Alec Mason and Wilson Stothers \cite{MAS3,MAS1} in the ring $R=\Int[i]$
of Gaussian integers.
\par
In this connection, it is very interesting to explicitly list
generators of the mixed commutator subgroups $[E(\Phi,R,A),E(\Phi,R,B)]$
{\it as subgroups\/}. From Theorem 5.1 we already know most of these
generators. These are $z_{\alpha}(\xi\zeta,\eta)$, where $\xi\in A$,
$\zeta\in B$, $\eta,\theta\in R$. But what are the remaining ones?

In fact, using the Chevalley commutator formula it is relatively
easy to show that $E(\Phi,R,A),E(\Phi,R,B)]$ is generated by its
intersections with the fundamental $SL_2$'s. Using somewhat more
detailed analysis the first and the fourth author established the
following result, initially for the case of $\GL(n,R)$, $n\ge 3$,
see \cite{RZ12} and then, jointly with the third author, for all
other cases, see \cite{yoga2,RNZ3}.

{\theorem Let $R$ be a commutative ring with $1$ and $A$, $B$
be two ideals of $R$. Then the mixed commutator subgroup
$\big[E(\Phi,R,A),E(\Phi,R,B)\big]$ is generated as a normal 
subgroup of $E(n,R)$ by the elements of the form

\begin{itemize}\itemsep=0pt
\item
$\big[x_{\alpha}(\xi),
{}^{x_{-\alpha}(\eta)}\!x_{\alpha}(\zeta)\big]$,
\item
$\big[x_{\alpha}(\xi),x_{-\alpha}(\zeta)\big]$,
\item
$x_{\alpha}(\xi\zeta)$,
\end{itemize}

\noindent
where $\alpha\in\Phi$, $\xi\in A$, $\zeta\in B$, $\eta\in R$.}

\bigskip

Another moderate technical effort allows to make it into a
natural candidate for the set of elementary generators of
$[E(\Phi,R,A),E(\Phi,R,B)]$.

{\theorem Let $R$ be a commutative ring with $1$ and $I$, $J$ be
two ideals of $R$. Then the mixed commutator subgroup
$\big[E(\Phi,R,A),E(\Phi,R,B)\big]$ is generated as a group by
the elements of the form

\begin{itemize}\itemsep=0pt
\item
$\big[z_{\alpha}(\xi,\eta),z_{\alpha}(\zeta,\theta)\big]$,
\item
$\big[z_{\alpha}(\xi,\eta),z_{-\alpha}(\zeta,\theta)\big]$,
\item
$z_{\alpha}(\xi\zeta,\eta)$,
\end{itemize}
\noindent
where $\alpha\in\Phi$, $\xi\in A$, $\zeta\in B$, $\eta,\theta\in R$.}


\subsection{Relative commutator width}

Now we are all set to address relative versions of the
main problem. The two following results were recently obtained
by the second author, with his method of universal localisation
\cite{step}, but they depend on the construction of generators
in Theorems 5.1 and 5.6. Mostly, the preceding results were
either published or prepublished in some form, and announced
at various conferences. These two theorems are stated here for
the first time.

{\theorem Let $R$ be a commutative ring with $1$ and let
$I\unlhd R$, be an ideal of $R$. Then there exists a natural
number $N=N(\Phi)$ depending on\/ $\Phi$ alone, such that
any commutator $[x,y]$, where
$$ x\in G(\Phi,R,I),\quad y\in E(\Phi,R)\qquad\text{or}\qquad
x\in G(\Phi,R),\quad y\in E(\Phi,R,I) $$
\noindent
is a product of not more that $N$ elementary generators
$z_{\alpha}(\xi,\zeta)$, $\alpha\in\Phi$, $\xi\in I$, $\zeta\in R$.}

{\theorem Let $R$ be a commutative ring with $1$ and let
$A,B\unlhd R$, be ideals of $R$. there exists a natural
number $N=N(\Phi)$ depending on\/ $\Phi$ alone, such
that any commutator
$$ [x,y],\qquad x\in G(\Phi,R,A),\quad y\in E(\Phi,R,B) $$
\noindent
is a product of not more that $N$ elementary generators
listed in Theorem $5.6$.}

\bigskip

Quite remarkably, the bound $N$ in these theorems does not depend
either on the ring $R$, or on the choice of the ideals $I,A,B$.
The proof of these theorems is not particularly long, but it
relies on a whole bunch of universal constructions and will be
published elsewhere. From the proof, it becomes apparent that 
similar results hold also in other such situations: for any 
other functorial generating set, for multiple relative 
commutators \cite{RZ12,RNZ3}, etc.


\section{Loose ends}

Let us mention some positive results on commutator width and
possible further generalisations.


\subsection{Some positive results}

There are some obvious bounds for the commutator width that follow
from unitriangular factorisations. For the $\SL(n,R)$ the following
result was observed by van der Kallen, Dennis and Vaserstein. The
proof in general was proposed by Nikolai Gordeev
and You Hong in 2005, but is still not published, as far as we know.

{\theorem Let $\rk(\Phi)\ge 2$. Then for any commutative ring $R$ an
element of $U(\Phi,R)$ is a product of not more than two commutators
in $E(\Phi,R)$.}

\bigskip

Combining the previous theorem with Theorem 3.4 we get the following
corollary.

{\corollary Let $\rk(\Phi)\ge 2$ and let $R$ be a ring such that
$\sr(R)=1$. Then the any element of $E(\Phi,R)$ is a product of
$\le 6$ commutators.}

\bigskip

This focuses attention on the following problem.

{\problem Find the shortest factorisation of $E(\Phi,R)$ of the
form
$$ E=UU^-UU^-\ldots U^{\pm}. $$}


Let us reproduce another result from the paper by Andrei Smolensky,
Sury and the third author \cite{VSS}. It is proven similarly to
Theorem 3.4, but uses Cooke---Weinberger \cite{CW} as induction base.
Observe that it depends on the Generalised Riemann's Hypothesis,
which is used to prove results in the style of Artin's conjecture
on primitive roots in arithmetic progressions. Lately, Maxim
Vsemirnov succeeded in improving bounds and in some cases eliminating
dependence on GRH. In particular, Cooke---Weinberger construct a
division chain of length 7 in the non totally imaginary case, the
observation that it can be improved to a division chain of length
5 is due to Vsemirnov \cite{vsemirnov}.
Again, in the form below, with $G(\Phi,{\mathcal O}_S)$ rather than
$E(\Phi,{\mathcal O}_S)$, it
relies on the almost positive solution of the congruence subgroup
problem \cite{bassmilnorserre,matsumoto}.

{\theorem Let $R={\mathcal O}_S$ be a Dedekind ring of arithmetic
type with infinite multiplicative group. Then under the Generalised
Riemann Hypothesis the simply connected Chevalley group
$G_{\sic}(\Phi,{\mathcal O}_S)$
admits unitriangular factorisation of length $9$,
$$ G_{\sic}(\Phi,{\mathcal O}_S)=UU^-UU^-UU^-UU^-U. $$
\noindent
In the case, where ${\mathcal O}_S$ has a real embedding, it
admits unitriangular factorisation of length $5$,
$$ G_{\sic}(\Phi,{\mathcal O}_S)=UU^-UU^-U. $$
\vskip -2truecm}

{\corollary Let $\rk(\Phi)\ge 2$ and let ${\mathcal O}_S$ be a
Dedekind ring of arithmetic type with infinite multiplicative group.
Then the any element of $G_{\sic}(\Phi,{\mathcal O}_S)$ is a product
of $\le 10$ commutators. In the case, where ${\mathcal O}_S$ has
a real embedding, this estimate can be improved to $\le 6$
commutators.}


\subsection{Conjectures concerning commutator width}

We believe that solution of the following two problems is now
at hand. In Section 2 we have already cited the works of
Frank Arlinghaus, Leonid Vaserstein, Ethel Wheland and You Hong
\cite{VW1,VW2,you94,AVY}, where this is essentially done for
{\it classical\/} groups, over rings subject to $\sr(R)=1$ or some
stronger stability conditions.

{\problem
Under assumption\/ $\sr(R)=1$ prove that any element of\/
$E_{\ad}(\Phi,R)$ is a product of\/ $\le 2$ commutators in\/
$G_{\ad}(\Phi,R)$.}

{\problem
Under assumption\/ $\sr(R)=1$ prove that any element of\/
$E(\Phi,R)$ is a product of\/ $\le 3$ commutators in\/
$E(\Phi,R)$.}

\bigskip

It may well be that under this assumption the commutator
width of $E(\Phi,R)$ is always $\le 2$, but so far we
were unable to control details concerning semisimple factors.
\par
It seems, that one can apply the same argument
to higher stable ranks. Solution of the following problem
would be a generalisation of \cite{DV88}, Theorem~4.

{\problem
If the stable rank\/ $\sr(R)$ of\/ $R$ is finite, and for
some $m\ge 2$ the elementary linear group\/
$E(m,R)$ has bounded word length with
respect to elementary generators, then for all\/ $\Phi$
of sufficiently large rank any element of\/
$E(\Phi,R)$ is a product of\/ $\le 4$ commutators in\/
$E(\Phi,R)$.}

{\problem
Let $R$ be a Dedekind ring of arithmetic type with infinite
multiplicative group. Prove that any element of\/
$E_{\ad}(\Phi,R)$ is a product of\/ $\le 3$ commutators in\/
$G_{\ad}(\Phi,R)$.}

\bigskip

Some of our colleagues expressed belief that any element
of $\SL(n,\Int)$, $n\ge 3$, is a product of $\le 2$ commutators.
However, for Dedekind rings with {\it finite\/} multiplicative
groups, such as $\Int$, at present we do not envisage any
{\it obvious\/} possibility to improve the generic bound $\le 4$
even for large values of $n$. Expressing elements of $\SL(n,\Int)$
as products of 2 commutators, if it can be done at all, should
require a lot of specific case by case analysis.


\subsection{The group $\SL(2,R)$: improved generators}

One could also mention the recent paper by Leonid Vaserstein
\cite{vaser10} which shows that for the group $\SL(2,R)$ it is natural
to consider bounded generation not in terms of the elementary
generators, but rather in terms of the generators of the
pre-stability kernel $\tilde E(2,R)$. In other words, one
should also consider matrices of the form $(e+xy)(e+yx)^{-1}$.

{\theorem The group\/ $\SL(2,\Int)$ admits polynomial parametrisation
of total degree\/ $\le 78$ with\/ $46$ parameters.}

\bigskip

The idea is remarkably simple. Namely, Vaserstein observes that
$\SL(2,\Int)$ coincides with the pre-stability kernel
$\tilde E(2,\Int)$. All generators of the group $\tilde E(2,\Int)$,
not just the elementary ones, admit polynomial parametrisation.
The additional generators require 5 parameters each.
\par
It only remains to verify that each element of $\SL(2,\Int)$
has a small length, with respect to this new set of generators.
A specific formula in \cite{vaser10} expresses an element of
$\SL(2,\Int)$ as a product of 26 elementary generators and
4 additional generators, which gives $26+4\cdot 5=46$ parameters
mentioned in the above theorem.


\subsection{Bounded generation and Kazhdan property}

The following result is due to Yehuda Shalom \cite{shalom99},
Theorem 3.4, see also \cite{shalom06,KN}.

{\theorem Let\/ $R$ be an $m$-generated commutative
ring,\/ $n\ge 3$. Assume that\/ $E(n,R)$ has bounded
width\/ $C$ in elementary generators. Then\/ $E(n,R)$ has
property\/ $T$. In an appropriate generating system\/
$S$ the Kazhdan constant is bounded from below
$$ {\mathcal K}(G,S)\ge \frac{1}{C\cdot22^{n+1}}. $$
\vskip -2truein}

{\problem Does the group\/ $\SL(n,\Int[x])$, $n\ge 3$,
has bounded width with respect to the set of elementary
generators?}

\bigskip

If this problem has positive solution, then by Suslin's
theorem and Shalom's theorem the groups\/ $\SL(n,\Int[x])$
have Kazhdan property\/ $T$. Thus,

{\problem Does the group\/ $\SL(n,\Int[x])$, $n\ge 3$,
have Kazhdan property\/ $T$?}

\bigskip

If this is the case, one can give a uniform bound of the
Kazhdan constant of the groups $\SL(n,{\mathcal O})$, for the
rings if algebraic integers. It is known that these group
have Kazhdan property, but the known estimates depend on the
discriminant of the ring $\mathcal O$.

{\problem Prove that the group\/ $\SL(n,\Rat[x])$ does not have
bounded width with respect to the elementary generators.}

\bigskip

It is natural to try to generalise results of Bernhard Liehl
\cite{liehl} to other Chevalley groups. The first of the
following problems was stated by Oleg Tavgen in \cite{tavgen90}.
As always, we assume that $\rk(\Phi)\ge 2$. Otherwise, Problem
6.10 is open for the group $\SL(2,{\mathcal O}_S)$, provided
that the multiplicative group ${\mathcal O}_S^*$ is infinite.

{\problem Prove that over a Dedekind ring of arithmetic type the
relative elementary groups\/ $E(\Phi,{\mathcal O}_S,I)$ have
bounded width with respect to the elementary generators
$z_{\alpha}(\xi,\zeta)$, with a bound that does not depend on $I$.}

{\problem Prove that over a Dedekind ring of arithmetic type
the mixed commutator subgroups\/ 
$[E(\Phi,{\mathcal O}_S,A),E(\Phi,{\mathcal O}_S,B)]$
have bounded width with respect to the elementary generators
constructed in Theorem $5.6$, with a bound that does not
depend on $A$ and $B$.}


\subsection{Not just commutators}

It is very challenging to understand, to which extent such
behaviour is typical for more general classes of group words.
There are a lot of recent results showing that the verbal length
of the finite simple groups is strikingly small \cite{Sh07,Sh09,LS,
LST,LOST12,GM}.
In fact, under some natural assumptions for large finite
simple groups this verbal length is 2. We do not expect
similar results to hold for rings other than the zero-dimensional
ones, and some arithmetic rings of dimension 1.
\par
Powers are a class of words in a certain sense opposite to
commutators. Alireza Abdollahi suggested that before passing
to more general words, we should first look at powers.

{\problem Establish finite width of powers in elementary
generators, or lack thereof.}

\bigskip

An answer -- in fact, {\it any\/} answer! -- to this problem
would be amazing. However, we would be less surprised if for
rings of dimension $\ge 2$ the verbal maps in $G(R)$ would
have very small images.


\section{Acknowledgements}

The authors thank Francesco Catino, Francesco de Giovanni and
Carlo Scoppola for an invitation to give a talk on commutator
width at the Conference in Porto Cesareo, which helped us to
focus thoughts in this direction. Also, we would like to thank
Nikolai Gordeev and You Hong for inspiring discussions of
positive results on commutator width and related problems,
Sury and Maxim Vsemirnov for discussion of arithmetic aspects,
Alireza Abdollahi for suggestion to look at powers, Anastasia
Stavrova amd Matthias Wendt for some very pertinent remarks
concerning localisation, and correcting some misprints in
the original version of our lemmas of conjugation calculus and
commutator calculus.

The second and the third authors started this research within the
framework of the RFFI/Indian Academy cooperation project
10-01-92651 ``Higher composition laws, algebraic $K$-theory and
algebraic groups'' (SPbGU--Tata Institute) and the RFFI/BRFFI
cooperation  project 10-01-90016 ``The structure of forms of
reductive groups, and behaviour of small unipotent elements in
representations of algebraic groups'' (SPbGU--Mathematics
Institute of the Belorussian Academy). Currently the work of the
second and the third authors is supported by the RFFI research
project 11-01-00756 (RGPU) and by the State Financed research
task 6.38.74.2011 at the Saint Petersburg State University
``Structure theory and geometry of algebraic groups and their
applications in representation theory and algebraic $K$-theory''.
The third author is also supported by the RFFI research project
12-01-00947 (POMI). The fourth author acknowledges the support
of NSFC grant 10971011.


\end{document}